\documentclass[12pt]{article}
\normalbaselines
\catcode `\@=11
\@addtoreset{equation}{section}

\def\section{\@startsection {section}{1}{\z@}{-3.5ex plus -1ex minus
     -.2ex}{2.3ex plus .2ex}{\normalsize\bf}}
\def\subsection{\@startsection{subsection}{2}{\z@}{-3.25ex plus -1ex minus
 -.2ex}{1.5ex plus .2ex}{\normalsize\bf}}
\def\thebibliography#1{\section*{References\markboth
  {REFERENCES}{REFERENCES}}\list
  {[\arabic{enumi}]}{\settowidth\labelwidth{[#1]}\leftmargin\labelwidth
  \advance\leftmargin\labelsep
  \usecounter{enumi}}
  \def\newblock{\hskip .11em plus .33em minus -.07em}
  \sloppy
  \sfcode`\.=1000\relax}
 
\catcode `\@=12

\def\ba{\begin{array}}
\def\ea{\end{array}}	
\def\be{\begin{equation}}
\def\ee{\end{equation}}
\def\bem{\begin{em}}
\def\eem{\end{em}}
\def\bpt{\begin{picture}}
\def\ept{\end{picture}}
\def\ot{\otimes}
\def\ld{\ldots}

\def\d{\delta}

\def\vp{\varphi}

\def\vp{\varphi}

\def\ra{\longrightarrow}

\def\ca{{\cal A}}

\def\ca{{\cal A}}

\def\cm{{\cal M}}

\def\cw{{\cal W}}

\def\lb{\langle}
\def\rb{\rangle}
\def\dn#1,{_{({#1})}}
\newfont{\numb}{msbm10}
\def\com{\mbox{\numb C}}

\def\pg#1,#2,#3,{\langle #1 | #2 \rangle^{#3}}
\def\otn#1,#2,{\left( #1 \otimes #2 \right) }

\def\ps#1,#2,{\Psi_{#1,#2}}
\def\id#1,{id_{#1}}
\def\ev#1,{ev_{#1}}
\def\tl#1,{\tilde{#1}}
\def\m#1,{m_{#1}}
\def\1n{^{(1)}}
\def\2n{^{(2)}}
\def\3{^{(3)}}
\def\ps#1,#2,{\Psi_{{#1}{,}{#2}}}
\def\et#1,{U^{\ot #1}}
\def\get#1,{{\cal U}^{\ast \ot #1}}

\begin{document}

\vspace*{2.5cm}
\noindent
{ \bf ON GENERALIZED STATISTICS AND INTERACTIONS}
\footnote{To be published in the Proceedings of the XVI 
Workshop on Geometric Methods in Physics, July 1-7, 1998, 
Bia{\l}owieza, Poland.
The work is partially sponsored by Polish Committee 
for Scientific Research (KBN) under Grant 2P03B130.12.}
\vspace{1.3cm}\\
\noindent
\hspace*{1in}
\begin{minipage}{13cm}
W{\l}adys{\l}aw Marcinek\vspace{0.3cm}\\
Institute of Theoretical Physics\\ 
University of Wroc{\l}aw, Poland\\
E-mail: wmar@ift.uni.wroc.pl
\end{minipage}

\vspace*{0.5cm}

\begin{abstract}
\noindent
The concept of exchange braid statistics is generalized.
The cross statistics is studied as a result of interaction. 
An algebraic model of a system of particles equipped with such 
statistics is described. The corresponding Fock space representation 
is also given.
\end{abstract}
%%%%%%%%%%%%%%%%%%%%%% section 1
\section{\hspace{-4mm}.\hspace{2mm}INTRODUCTION}
%\subsection{\hspace{-5mm}.\hspace{2mm}First Subsection of the Current
%Section (main words capitalized)}
It is well known that the standard approach to statistics
of identical particles is based on the notion of the usual 
symmetric group $S_n$. This group describes the interchange 
process of indistinguishable particles. The transposition
$(1, 2) \ra (2, 1)$ corresponding for the exchanging of two 
identical particles is represented by the map
\be
\ba{c}
\tau : x^1 \ot x^2 \ra \pm \, x^2 \ot x^1. 
\label{tra}
\ea
\ee
Every transposition yields a phase factor equal to $\pm 1$ 
($+ 1$ for bosons and $-1$ for fermions). If we replace the factor 
$\pm 1$ by a complex parameter $q$, then we obtain the most simple 
generalization of the usual concept of statistics, namely the 
well--known $q$--statistics \cite{owg,gre,moh,mepe,bs2}.
The corresponding particles are said to be quons. If 
$q := \exp(i \vp)$, where $0 \leq \vp < 2\pi$ is the 
so--called statistics parameter, then the corresponding 
$q$--statistics is determined by the value of $\vp$. 
Observe that for $\vp = 0$ we have bosons,
and for $\vp = \pm \pi$ -- fermions. For arbitrary 
$\vp \in [0, 2\pi)$ we obtain anyons \cite{wil,wiz}. 

There is an interesting concept of an exotic statistics in a low 
dimensional space based on the notion of the braid group $B_n$ 
\cite{Wu,I}. In this concept the configuration space for the 
system of $n$--identical particles moving on a manifold $\cm$ 
is given by the following relation
$$
Q_{n}(\cm) = {\left(\cm^{\times n} - D\right) }/{S_{n}},
$$ 
where $D$ is the subset of the Cartesian product $\cm^{\times n}$ 
on which two or more particles occupy the same position. 
The group $\pi_{1}\left( Q_{n}(M)\right)\equiv B_{n} (M)$
is just the $n$--string braid group on $\cm$. Note that the 
statistics of a system of particles is determined by the group 
$\Sigma_n$ \cite{Wu,I}. This group is a subgroup of the braid 
group $B_n(\cm)$ corresponding for the interchange process of 
two arbitrary indistinguishable particles. It is an extension 
of the symmetric group $S_n$. The mathematical formalism related 
to the braid group has been developed intensively by Majid, see
\cite{SM,Maj,SMa,bm,sma,qm,ssm} for example. An algebraic 
formalism for a particle system with generalized statistics 
has been considered by the author \cite{WM3,WM4,WM6,WM8,wm7}.
It is interesting that in this algebraic approach all 
commutation relations for particles equipped with an 
arbitrary statistics can be described as a representation 
of the so--called quantum Weyl algebra $\cw$ (or Wick algebra)
\cite{jswe,RM,m10,ral,bor}. In this attempt the creation and 
annihilation operators act on an algebra $\ca$. The creation operators 
act as the multiplication in this algebra and the annihilation ones act 
as a noncommutative contraction (noncommutative partial derivatives). 
The algebra $\ca$ play the role of noncommutative Fock space. The 
application for particles in singular magnetic field has been given 
by the author \cite{mco,top,sin}. Note that similar approach has been 
also considered by others authors \cite{sci,twy,mep,mphi,mira}. 
An interesting concept related to generalized statistics has been 
also given in \cite{fios,melme}.

In this paper we are going to study of a system of charged particles  
moving under influence of an intermediate quantum field. Our 
fundamental assumption is that every charged particle is transform 
under interaction into a system consisting a charge and  quanta of 
the field. Such system behaves like free particles moving in certain 
effective space.  It is showed that the system is endowed with the 
so--called cross statistics. This statistics is a result of 
inetractions. For the description of such statistics 
we develop the algebraic model of a system with generalized 
statistics studied previously by the author \cite{castat,qstat}.
%%%%%%%%%%%%%%%%%%%%%%% section 2
\section{\hspace{-4mm}.\hspace{2mm}FUNDAMENTAL ASSUMPTIONS}
We are going here to study a system of charged particles with
certain dynamical interaction. It is natural to expect that some 
new and specific quantum states of the system have appear as a 
result of interaction. We would like to describe all such states. 
In order to do this we assume that the interaction can be described 
by an intermediate quantum field. 
Our fundamental assumption is that every charged particle is 
transform under interaction into a system consisting a charge 
and $N$--species of quanta of the field. A system which contains 
a charge and certain number of quanta as a result of interaction 
with the quantum field is said to be a {\it dressed particle}
\cite{gso}. Next we assume that every dressed particle is a composite
object equipped with an internal structure. Obviously the structure
is determined by the interaction with the quantum field. We describe 
the structure of a dressed particle as a nonlocal system which
contains $n$ centers. Such centers behave like free particles
moving on certain effective space. Every center is equipped with 
ability for absorption and emission of quanta of the intermediate 
field. A center dressed with a single quantum of the field is said
to be a {\it quasiparticle}. Quasiparticles represent elementary 
excited quantum states of the given system and they are are described 
by as a finite set of elements
\be 
Q := \{x^i : i = 1,\ldots, N<\infty\}
\ee
which form a basis for a linear space $E$ over a field of complex 
numbers $\com$. A center which is an empty place for a quantum is 
called a {\it quasihole}. Quasiholes represent conjugate states 
and they are corresponding to the basis 
\be
Q^{\ast} := \{x^{\ast i} : i = N, N-1,\ldots , 1\}.
\ee
for the complex conjugate space $E^{\ast}$. The pairing 
$g_E : E^{\ast} \ot E \ra \com$ and the corresponding
scalar product is given by
\be
g_E (x^{\ast i}\ot x^j) = \lb x^{i}|x^j \rb := \d^{ij}.
\ee
A center which contains any quantum is said to be {\it neutral}.
It represents the ground state $|0> = {\bf 1}$ of the system.
A neutral center can be transform into a quasiparticle or a
quasihole by an absorption or emission process of single quantum,
respectively. 
Quasiparticles and quasiholes as components of certain dressed
particle have also their own statistics. It is interesting that  
there is a statistics of new kind, namely a {\it cross statistics}. 
This statistics is determined by an exchange process of quasiholes
and quasiparticles. Note that the exchange is not a real process 
but an effect of interaction. Such exchange means annihilation
of a quasiparticle on certain place and simultaneous creation
of quasihole on an another place. This statistics is described 
by an operator $T$ called an {\it elementary cross} or {\it twist}.
This operator is linear, invertible and Hermitian. It is given by
its matrix elements
$T : E^{\ast}\ot E\ra E \ot E^{\ast}$
\be
\ba{c}
T(x^{\ast i}\ot x^j) = \Sigma \ T^{ij}_{kl} x^k \ot x^{\ast l}.
\label{cross}
\ea
\ee
The usual exchange statistics of quasiparticles is described a
linear $B$ satisfying the standard braid relations
\be
\ba{c}
B\1n B\2n B\1n = B\2n B\1n B\2n ,
\label{bra}
\ea
\ee
where $B\1n := B \ot id$ and $B\2n := id \ot B$. The exchange
process determined by the operator $B$ is a real process. Such
exchange process is possible if the dimension of the effective
space is equal or great than two. Hence in this case we need 
two operators $T$ and $B$ for the description of our system with 
generalized statistics. These operators are not arbitrary. 
They must satisfy the following consistency conditions
\be
\ba{c}
B^{(1)}T\2nT\1n = T\2nT\1n B^{(2)},\\
(id_{E \ot E} + \tilde{T})(id_{E \ot E} - B) = 0,
\label{cod}
\ea
\ee
where the operator $\tilde{T} : E\ot E\ra E\ot E$ is given 
by its matrix elements 
\be
(\tilde{T})^{ij}_{kl} = T^{ki}_{lj}.
\ee
We need a solution of these conditions for the construction
of an example of a system with generalized statistics.
Note that the general solution for these conditions is 
not known. Hence we must restrict our attention to some 
particular cases only. One can use solutions used in 
noncommutative differential calculi in order to give some 
examples \cite{quo}. In the one--dimensional case
there is no place for such exchange process. Hence in this case
only the cross statistics is possible, the exchange braid 
statistics not exists. 
%%%%%%%%%%%%%%%%%%%%%%%%%%%%%%%%%%%%%%%%%%
\section{\hspace{-4mm}.\hspace{2mm}HERMITIAN WICK ALGEBRAS}
Let us consider a pair of unital and associative algebras $\ca$ 
and $\ca^{\ast}$. We assume that they are conjugated. This means 
that there is an antilinear and involutive isomorphism 
$(-)^{\ast} : \ca\ra\ca^{\ast}$ and we have the following relations
\be
m_{\ca^{\ast}}(b^{\ast} \ot a^{\ast}) 
= (m_{\ca}(a \ot b))^{\ast},
\quad (a^{\ast})^{\ast} = a,
\ee
where $a, b \in \ca$ and $a^{\ast}, b^{\ast}$ are their images 
under the isomorhism $(-)^{\ast}$. Both algebras $\ca$ and 
$\ca^{\ast}$ are graded
\be
\ba{cc}
\ca := \bigoplus\limits_{n} \ \ca^n ,&
\ca^{\ast} := \bigoplus\limits_{n} \ \ca^{\ast n} .
\ea
\ee
A linear mapping $\Psi : \ca^{\ast}\ot\ca\ra\ca\ot\ca^{\ast}$ 
such that we have the following relations 
\be
\ba{l}
\Psi \circ (id_{\ca^{\ast}} \ot m_{\ca}) 
= (m_{\ca} \ot id_{\ca^{\ast}}) 
\circ (id_{\ca} \ot \Psi) \circ (\Psi \ot id_{\ca}),\\
\Psi \circ (m_{\ca^{\ast}} \ot id_{\ca}) 
= (id_{\ca} \ot m_{\ca^{\ast}}) 
\circ (\Psi \ot id_{\ca^{\ast}}) 
\circ (id_{\ca^{\ast}} \ot \Psi)\\
(\Psi(b^{\ast}\ot a))^{\ast} = \Psi (a^{\ast} \ot b)
\label{twc}
\ea
\ee
is said to be a cross symmetry or $\ast$--twist \cite{bma}. 
We use here the notation
\be
\Psi(b^{\ast}\ot a) = \Sigma a_{(1)}\ot b^{\ast}_{(2)}
\ee
for $a\in\ca, b^{\ast}\in\ca^{\ast}$.

The tensor product $\ca\ot\ca^{\ast}$ of algebras $\ca$ and
$\ca^{\ast}$ equipped with the multiplication 
\be
\ba{c}
m_{\Psi} := (m_{\ca} \ot m_{\ca^{\ast}}) \circ 
(id_{\ca} \ot \Psi \ot id_{\ca^{\ast}})
\label{mul}
\ea
\ee
is an associative algebra called a Hermitian Wick algebra 
\cite{jswe,bma} and it is denoted by
$\cw = \cw_{\Psi}(\ca) = \ca \ot_{\Psi} \ca^{\ast}$. 
This means that the Hermitian Wick algebra $\cw$ is the
tensor cross product of algebras $\ca$ and $\ca^{\ast}$
with respect to the cross symmetry $\Psi$ \cite{bma}.
Let $H$ be a linear space. We denote by $L(H)$ the algebra of 
linear operators acting on $H$. One can prove \cite{bma}
that we have the

{\bf Theorem:} 
Let $\cw \equiv \ca \ot_{\Psi} \ca^{\ast}$  be a Hermitian Wick 
algebra. If $\pi_{\ca} : \ca \ra L(H)$ is a 
representation of the algebra $\ca$, such that we have the relation
\be
\ba{c}
(\pi_{\ca}(b))^{\ast} \pi_{\ca}(a) = \Sigma
\pi_{\ca}(a_{(1)}) \pi_{\ca^{\ast}}(b^{\ast}_{(2)}),\\
\pi_{\ca^{\ast}}(a^{\ast}) := (\pi_{\ca}(a))^{\ast},
\label{wre}
\ea
\ee
then there is a representation $\pi_{\cw} : \cw \ra L(H)$ of the 
algebra $\cw$. 
We use the following notation
\be
\pi_{\ca}(x^i) \equiv a_{x^i}^+ , \quad
\pi_{\ca^{\ast}}(x^{\ast i}) \equiv a_{x^{\ast i}} .
\ee
The relations (\ref{wre}) are said to be {\it commutation relations} 
if there is a positive definite scalar product on $H$ such that 
operators $a_{x^i}^+$ are adjoint to $a_{x^{\ast i}}$ and vice versa.
Let us consider a Hermitian Wick algebra $\cw$ corresponding for a 
system with generalized statistics. For the construction of such 
algebra we need a pair of algebras $\ca$, $\ca^{\ast}$ and a cross 
symmetry $\Psi$. It is natural to assume that these algebras have
$E$ and $E^{\ast}$ as generating spaces, respectively, and there is
the following condition for the cross symmetry 
\be
\Psi|_{E^{\ast}\ot E} = T + g_E .
\ee
%%%%%%%%%%%%%%%%%%%%%%%%%%%%%%%%%%%%%%%%%%%%%%%%%%%%%%%%%%%%%
\section{\hspace{-4mm}.\hspace{2mm}FOCK SPACE REPRESENTATION}
Let us consider the Fock space representation of the algebra
$\cw$ corresponding for a system with generalized statistics. 
For the ground state and annihilation operators we assume that
\be
\langle 0|0 \rangle = 0, \quad a_{s^{\ast}} |0\rangle = 0 
\quad \mbox{for} \quad s^{\ast} \in \ca^{\ast}.
\ee
In this case the representation act on the algebra $\ca$.
Creation operators are defined as a multiplication in the 
algebra $\ca$
\be
a^+_{s} t := m_{\ca}(s \ot t), \quad 
\mbox{for} \quad s, t \in \ca .
\ee
The proper definition of the action of annihilation operators 
on the whole algebra $\ca$ is a problem. 

If the action of 
annihilation operators are given in such a way that there 
is unique, nondegenerate, positive definite scalar product
on $\ca$, creation operators are adjoint to annihilation ones 
and vice versa, then we say that we have a well--defined 
system with generalized statistics in the Fock representation
\cite{qstat}. 

Let us consider some examples for such systems. Assume that
quasiparticles and quasiholes are moving on one dimensional
effective space. In this case the algebra of states $\ca$ is 
the full tensor algebra $TE$ over the space $E$, and the conjugate
algebra $\ca^{\ast}$ is identical with the tensor algebra 
$TE^{\ast}$. If $T \equiv 0$ then  we obtain the most simple 
example of well--defined system with generalized statistics. 
The corresponding statistics is the so--called infinite 
(Bolzman) statistics \cite{owg,gre,qstat}. If 
$T : E^{\ast}\ot E\ra E \ot E^{\ast}$ is an arbitrary 
nontrivial cross operator, then there is the cross symmetry 
$\Psi^T : TE^{\ast}\ot TE\ra TE\ot E^{\ast}$. 
It is defined by a set of mappings 
$\Psi_{k,l}:E^{\ast\ot k}\ot E^{\ast\ot k}
\ra E^{\ot l}\ot E^{\ast\ot k}$,
where $\Psi_{1,1} \equiv R := T + g_E$, and
\be
\ba{l}
\Psi_{1,l} := R^{(l)}_l \circ \ld \circ R^{(1)}_l,\\
\Psi_{k,l} := (\Psi_{1,l})^{(1)}  \circ \ldots \circ 
(\Psi_{1,l})^{(k)},
\label{up}
\ea
\ee
here $R^{(i)}_l : E^{(i)}_l \ra E^{(i+1)}_l$,
$E^{(i)}_l := E \ot \ldots \ot E^{\ast} \ot E \ot \ldots \ot E$
($l+1$-factors, $E^{\ast}$ on the i-th place, $i \leq l$) 
is given by the relation
$$
R^{(i)}_l := \underbrace{
id_{E} \ot \ldots \ot R \ot \ldots \ot id_{E}
}_{l\;\;times},
$$
where $R$ is on the i-th place, $(\Psi_{1,l})^{(i)}$ is defined in 
similar way like $R^{(i)}$. The commutation relations (\ref{wre}) 
can be given here in the following form
\be
a_{x^{\ast i}} a_{x^j}^+ - T^{ij}_{kl} \ 
a_{x^k}^+ a_{x^{\ast l}} = \delta^{ij}{\bf 1}.
\ee
If the linear operator $\tilde{T} : E\ot E\ra E\ot E$ 
with the following matrix elements
\be
(\tilde{T})^{ij}_{kl} = T^{ki}_{lj}.
\ee
is bounded, we have the following Yang-Baxter equation on 
$E\ot E \ot E$ 
\be
(\tl T, \ot id_E )\circ (id_E \ot \tl T, )\circ 
(\tl T, \ot id_E ) = (id_E \ot \tl T, ) \circ 
(\tl T, \ot id_E )\circ (id_E \ot \tl T, ),
\ee
and $||\tl T,|| \leq 1$, then according to Bo$\dot{z}$ejko and 
Speicher \cite{bs2} there is a positive definite scalar product.
Note that the existence of a nontrivial kernel of the operator 
$P_2 \equiv id_{E\ot E} + \tilde{T}$ is essential 
for the nondegeneracy of the scalar product \cite{jswe}.
One can see that if this kernel is trivial, then we
obtain well--defined system with generalized statistics
\cite{RM,ral}.

If the dimension of the effective space is great than one
and the kernel of $P_2$  is nontrivial, then the scalar product 
is degenerate. Hence we must remove this degeneracy by factoring 
the mentioned above scalar product by the kernel. In this case 
we have $\ca := TE/I, \quad \ca^{\ast} := TE^{\ast}/I^{\ast}$,
where $I := gen\{id_{E\ot E} - B\}$ is an ideal in $TE$ and
$B: E \ot E \ra E \ot E$ is a linear and invertible operator 
satisfying the braid relation (\ref{bra}) and the consistency 
conditions (\ref{cod}), $I^{\ast}$ is the corresponding conjugated 
ideal in $TE^{\ast}$. One can see that there is the cross symmetry 
and the action of annihilation operators can be defined
in such a way that we obtain the well--defined system with 
the usual braid statistics \cite{RM,ral}. We have here the 
following commutation relations 
\be
\ba{l}
a_{x^{\ast i}} a_{x^j}^+ - T^{ij}_{kl} \ 
a_{x^k}^+ a_{x^{\ast l}} = \delta^{ij}{\bf 1}\\
a_{x^{\ast i}} a_{x^{\ast j}} - B_{ij}^{kl} \ 
a_{x^{\ast l}} a_{x^{\ast k}} = 0, \\
a^+_{x^i} a^+_{x^j} - B^{ij}_{kl} \ 
a^+_{x^k} a^+_{x^l} = 0.
\ea
\ee
Observe that for $T \equiv B \equiv \tau$, where $\tau$ represents 
the transposition (\ref{tra}) we obtain the usual canonical
(anti--)commutation relations for bosons or fermions.

\section*{Acknowledgments}  
The author would like to thank A. Borowiec for the discussion
and any other help.

\end{document}